\documentclass[preprint, 12pt, draft]{elsarticle}

\usepackage[centertags]{amsmath}
\usepackage{multirow,newlfont,amssymb,amsfonts,amsthm,booktabs,longtable,xcolor}

\DeclareMathOperator{\Span}{span}

\newdefinition{defn}{Definition}[section]

\theoremstyle{remark}
\newtheorem{rem}[defn]{Remark}

\makeatletter
\def\ps@pprintTitle{%
  \let\@oddhead\@empty
  \let\@evenhead\@empty
  \def\@oddfoot{\reset@font\hfil\thepage\hfil}
  \let\@evenfoot\@oddfoot
}
\makeatother

\allowdisplaybreaks[1] 

\begin{document} 


\begin{frontmatter}


\title{Group Classification of a generalization of the  Heath Equation}


\author[imecc]{Y.~Bozhkov}
\ead{bozhkov@ime.unicamp.br}

\author[imecc]{S.~Dimas\corref{cor}}
\ead{spawn@math.upatras.gr}

\cortext[cor]{Corresponding author}


\address[imecc]{Instituto de Matem\'atica, Estat\'istica e Computa\c{c}\~ao Cient\'ifica - IMECC\\ 
Universidade Estadual de Campinas - UNICAMP, Rua S\'ergio Buarque de Holanda, $651$\\
$13083$-$859$ - Campinas - SP, Brasil}


\begin{abstract}
	The complete group classification of a generalization of the Heath model is carried out by connecting it to the heat equation with nonlinear source.  Examples of invariant solutions are given under the terminal and the barrier option condition.
  \end{abstract}


\end{frontmatter}

\section{Introduction}

The area of Financial Mathematics was founded by the pioneering works of Samuelson, Modigliani, Merton, Black and Scholes with the formal introduction of the concepts of stochastic calculus into Econometrics. Through the modeling process of incorporating the randomness into a deterministic model a set of evolutionary PDEs is obtained. Amongst them are:
\begin{itemize}[-]
	\item The Black--Scholes--Merton Equation \cite{BlaScho73},
		\begin{equation}\label{eq:BSM}
			u_t+\frac{1}{2}\sigma^2x^2u_{xx}+rxu_x-ru=0,
		\end{equation}
	\item The Longstaff Equation \cite{Lo89},
		\begin{equation}\label{eq:Longstaff}
			u_t-\left(\frac{1}{4}\sigma^2-\kappa\sqrt{x}-2\lambda x\right)u_x-\frac{1}{2}\sigma^2x u_{xx}+xu=0,
		\end{equation}
	\item The Vasicek Equation \cite{Va77},
		\begin{equation}\label{eq:Vasicek}
			u_t+\frac{1}{2}\sigma^2u_{xx}-(\kappa-\lambda x)u_x-xu=0,
		\end{equation}	
	\item The Cox--Ingersoll--Ross Equation \cite{CoIngeRo85},
		\begin{equation}\label{eq:CIR}
			u_t+\frac{1}{2}\sigma^2xu_{xx}-(\kappa-\lambda x)u_{x}-x u=0,
		\end{equation}
	\item A Hamilton--Jacobi--Bellman type Equation \cite{HeaPlaSchwei2k1},
		\begin{equation}\label{eq:Heath}
			u_t+\alpha u_x+\frac{1}{2}b^2u_{xx}-\frac{1}{2}u^2_x+\nu(x)=0,
		\end{equation}
\end{itemize} 
where $u=u(x,t)$. Eqs.~\eqref{eq:BSM}--\eqref{eq:CIR} have similar structure, in fact, they are connected with the Heat Equation, 
\begin{equation}\label{eq:heat}
	u_t=u_{xx},
\end{equation}
via an invertible transformation as first shown for Eq.~\eqref{eq:BSM} in \cite{GazIbr98} and later for the rest of them in \cite{DiAndTsLe2k9}. Actually, Eqs.~\eqref{eq:Longstaff}--\eqref{eq:CIR} are encompassed in a general bond-pricing equation  presented in \cite{SiLeaHa2k8}. Furthermore, it is shown that Eqs.~\eqref{eq:BSM}--\eqref{eq:CIR} can be contained in a nonlinear generalization of this bond-pricing equation \cite{BoDi2k13x}.

The last of the above equations,\eqref{eq:Heath}, was first presented by Heath et al in  \cite{HeaPlaSchwei2k1} as an equation for mean variance hedging. Even though at first glance it looks to be different in structure than the rest, it can be linearized, and actually, for some specific cases of $\nu(x)$ can  be linked with \eqref{eq:heat} as well \cite{NaiAndrioLea2k5, DiAndTsLe2k9}.

For all the aforementioned connections of the Eqs.~\eqref{eq:BSM}--\eqref{eq:CIR} with the heat equation the use of symmetries  was prominent. The symmetry analysis of differential equations is a method first developed in the $19^\text{th}$ century by Sophus Lie. One of the main benefits of this method is that by following a completely algorithmic procedure one is able to determine the symmetries of a differential equation or systems of differential equations. \textit{Grosso modo}, the symmetries of a differential equation transform solutions of the equation to other solutions. The Lie point symmetries comprise a structural property of the equation --- in essence they are the  DNA of an equation. The knowledge of the symmetries of an equation enables one to utilize them for a variety of purposes, from obtaining analytical solutions and reducing its order to finding of integrating factors and conservation laws. In fact, many, if not all, of the different empirical methods for solving ordinary differential equations (ODEs) we have learned from standard courses at the undergraduate level emerge from a symmetry. For instance, having at our disposal a Lie point symmetry of a first order ODE, we can immediately get explicitly an integrating factor. Furthermore, even the knowledge of a trivial solution of the equation can be used for creating nontrivial solutions by using the equation's symmetries. And all these are due to the rich underlining algebraic structure of Lie groups and Lie algebras with which we give flesh to the symmetries of a differential equation. 

Another important characteristic of the symmetry method is that in some situations the symmetries of an equation may indicate that it can be transformed to a linear equation, to be linearized. In addition, its symmetries provide the means to construct the needed transformation. A strong indication for that is the existence of an infinite dimensional Lie algebra, \cite{BluKu89}.  

Moreover,  a valuable tool when considering classes of equations containing arbitrary constants and functions is the use of the \emph{equivalence} or \emph{admissible transformations}, \cite{Ovsi82,Ibra2k9,PoEshra2k5}. Equivalence transformations of a class of differential equations are point transformations that keep this class invariant, in other words they map an equation from this class  to another member of the same class. In the recent years equivalence transformations have found much application either as a stand alone analytic  tool for the group classification of differential equations, \cite{RoTo99, BoDi2k13b,BoDi2k13x}, or at the core of the \emph{enhanced group analysis},  \cite{CheSeRa2k8,CaBiPo2k11,IvaPoSo2k10,VaPoSo2k9}. 

In the present work the following generalization of Eq.~\eqref{eq:Heath},
\begin{equation}\label{eq:NonLinHeath}
	u_t=a u_x+\frac{1}{2}u^2_x-\frac{1}{2}b^2u_{xx}+f(x,u),\ b,f_u\ne0,
\end{equation}
 is studied under the prism of symmetry analysis. It is evident that in general Eq.~\eqref{eq:NonLinHeath} will not be linearizable as Eq.~\eqref{eq:Heath}. 
 
 Furthermore, two distinct conditions will be considered along with \eqref{eq:NonLinHeath}. The terminal condition  
\begin{equation}\label{cond:terminal}
	u(x,T)=1, 
\end{equation}
where $T$ is the terminal time, and the barrier option condition
\begin{subequations}\label{cond:barrier}
	\begin{align}
		&u(H(t),t)=R(t),\label{cond:barriera}\\ 
		&u(x,T)=\max(x-K,0),	\label{cond:barrierb}
	\end{align}
\end{subequations}
where the barrier option $u(x,t)$ satisfies Eq.~\eqref{eq:NonLinHeath} for $x>H(t),\ t<T$.\footnote{It is evident that the terminal condition \eqref{cond:terminal} is included as a special situation when the asset price exceeds a value defined by the barrier function $H$ on time $T$.} The constant  $T$ again is the terminal time where the barrier option is exercised and $K$ is the strike price.

The  first condition (\ref{cond:barriera}) describes the evolution of standard or ``vanilla"  products; the price of a zero-coupon bond (or of a financial option), $u(x,t)$, which expires when  $t=T$ \cite{BlaScho72,BlaScho73,Me74,Ugur2k8}. 

The second condition (\ref{cond:barrierb}) takes into account the possibility of an unacceptable change in the price of the underlying option. It is considered an exotic type of option possessing features that makes it more complex than  the ``vanilla" option, \cite{ha2k11, Kwo2k8, HaSoLea2k13}:  The underlying idea is that now a barrier $H(t)$ exists and when the asset price $x$ crosses it, the barrier option $u(x,t)$  becomes extinguished or  comes into existence. Those two types are also known as  \emph{down-and-out} and \emph{down-and-in}, respectively. Often a rebate, $R(t)$, is paid if the option is extinguished.  In what follows we shall consider the down-and-out type. A common assumption for the barrier function $H$ is to have the exponential form
\begin{equation}\label{H0}
H(\tau) = \beta Ke^{-\alpha \tau},
\end{equation}
where $\alpha\ge0$, $0\le \beta\le1$ and $\tau=t-T$ \cite[p. 187]{Kwo2k8}.

The main purpose of this paper is to carry out the complete group classification of Eq.~\eqref{eq:NonLinHeath} and  construct analytical solutions satisfying one of the two conditions. Recall that to perform a complete group classification of a differential equation (or a system of differential equations) involving arbitrary functions or/and parameters, means to find the Lie point symmetry group $\mathcal{G}$ for the most general case, and then to find specific forms of the differential equation for which $\mathcal{G}$ can be enlarged, \cite[p. 178]{Olver2k}. Quite often, there is good physical or geometrical motivation to study such cases. Indeed, the larger the group the bigger the number of symmetries that will ``survive" after imposing the conditions and hence the greater the chance of constructing an analytical solution corresponding to a more realistic financial model --- in our case --- following more closely the complex evolution of the financial option or bond. This is particularly true for the case of the barrier option as recent results have shown \cite{HaSoLea2k13, BoDi2k13b}.
 
As previously mentioned the key analytical tool used in this work is the symmetry analysis of Eq.~\eqref{eq:NonLinHeath}. One of the advantages of this approach, as already emphasized, is that it provides a well defined algorithmic procedure which essentially enables one to find the involved linearizing transformations, conservation laws, invariant solutions, etc. On the other hand, the calculations involved are usually very difficult and extensive even for the simplest equations.  Thus, it may become very tedious and error prone. For this reason the real progress in this area occurred in the last few decades with the advances in computer technology and the development of computer algebra systems like {\it Mathematica, Maple, Reduce}, etc.  Based on these systems, a handful of symbolic packages for determining the symmetries of differential equations exists, \cite{Head93,Nucci1,Nucci2,Baumann2k}. Such a symbolic package,  based on {\it Mathematica}, \cite{Wolfram2k10}, has been developed by S. Dimas. The package, named SYM, \cite{Dimas2k8,DiTs2k5a,DiTs2k6}, was developed from the ground up using the symbolic manipulation power of  {\it Mathematica} and the artificial intelligence capabilities  which it offers.  It was extensively used for all the results in the present paper, both for the interactive manipulation of the found symmetries and for the classification of the equations employing the symbolic tools provided by it. 

This paper is organized as follows. In section 2 the basic concepts of the Lie point symmetry approach to differential equations used in the paper are presented. In section 3 we show how Eq.~\eqref{eq:NonLinHeath} can be connected with the heat equation with nonlinear source and hence reduce the problem of the complete group classification to it. In section 4  we give various examples of how the classification can be used in order to obtain similarity solutions under the specific boundary problems studied, the ``vanilla" option and the barrier option.  Finally, in section 5 the obtained results in this work are discussed.

\section{Preliminaries}

In this section we expose some notions of the modern group analysis that will be encountered in the main sections of the article suitably adapted to the article's needs. For a full treatise of the subject we direct the interested reader to the classical texts  \cite{Olver2k,BluKu89,Ibra85b,Ovsi82,Hydon2k,Ste90}.

A Lie point symmetry\footnote{Henceforth referred simply as symmetry.} of \eqref{eq:NonLinHeath} is a differential operator, named \emph{infinitesimal generator},
 \begin{equation*}\label{prelim:symmetryOperator}
  	\mathfrak{X}={\xi }^1(x,t,u)\frac{\partial }{\partial x}+{\xi }^2(x,t,u)\frac{\partial }{\partial t}+ \eta(x,t,u)  \frac{\partial }{\partial u},
\end{equation*}
that satisfy the condition
\begin{equation}\label{prelim:linearizedsymmetrycondition}
	\left.\mathfrak{X}^{(2)}\left[a u_x+\frac{1}{2}u^2_x-\frac{1}{2}b^2u_{xx}+f(x,u)-u_t\right]\right\rvert_\eqref{eq:NonLinHeath} \equiv 0,
\end{equation}
where $\mathfrak{X}^{(2)}$ is the suitable prolongation of the differential operator up to order two. The set of all the symmetries admitted by a differential equation constitutes a Lie algebra.

Having the symmetries there is a wealth of things that can be done with. In the present paper, we use them to obtain \emph{invariant} or \textit{similarity solutions} of the Eq.~\eqref{eq:NonLinHeath}. By invariant solutions we mean solutions of \eqref{eq:NonLinHeath} that are invariant under one of the found symmetries $\mathfrak{X}$, that is, 
\begin{equation}\label{prelim:invariantsurfacecond}
	\left.\mathfrak{X}[u-\varphi (x,t)]\right\rvert_{u=\varphi (x,t)}\equiv0.
\end{equation}
The Eq.~\eqref{prelim:invariantsurfacecond} is a linear PDE  called \textit{invariant surface condition} and by solving it we obtain a way to reduce the order of Eq.~\eqref{eq:NonLinHeath}. Similarly, when we look for a similarity solution of Eq.~\eqref{eq:NonLinHeath} along with a initial/boundary condition we have to choose the subalgebra  leaving also invariant that condition and its boundary:
\begin{equation}\label{prelim:terminala}
	\left.\mathfrak{X}(t-T)\right\rvert_{t=T}\equiv0
\end{equation}
and
\begin{equation}\label{prelim:terminalb}
	\left.\mathfrak{X}(u-1)\right\rvert_{t=T,\,u=1}\equiv0
\end{equation}
for the boundary condition \eqref{cond:terminal}. And
\begin{equation}\label{prelim:barriera}
	\left.\mathfrak{X}(x-H(t))\right\rvert_{x=H(t)}\equiv0,
\end{equation}
\begin{equation}\label{prelim:barrierb}
	\left.\mathfrak{X}(u-R(t))\right\rvert_{x=H(t),\, u=R(t)}\equiv0
\end{equation}
for the boundary condition \eqref{cond:barrier}.

Another useful notion is that of \emph{additional equivalence transformation}. An additional equivalence transformation is a point transformation that connects inequivalent --- from the point of view of equivalence transformations --- subclasses of a broader class of differential equations, which contains the class of differential equations under consideration, that possess the same Lie algebra of symmetries. The knowledge of such transformations greatly facilitates the classification.

\section{Group classification}

In this section we proceed with the group classification of Eq.~\eqref{eq:NonLinHeath}.  First, we observe that the point transformation
\begin{equation}\label{main:PointTransformations}
	\tau =-\frac{b^2}{2}t,\ \phi = \exp(-\frac{1}{b^2}(a x+u)),
\end{equation}
turns the class \eqref{eq:NonLinHeath} into the heat equation with nonlinear source,
\begin{equation}\label{main:heatNonLinearSource}
	\phi_\tau=\phi_{x x} + \hat f(x, \phi),\ \hat f_{\phi\phi}\ne0,
\end{equation}
where $\hat f = \frac{1}{b^4}\exp(-\frac{1}{b^2}(a x+u))(2f -a^2)$. Hence the problem of the complete group classification of the class \eqref{eq:NonLinHeath}  is reduced to the group classification of the heat equation with nonlinear source, \eqref{main:heatNonLinearSource}.
\begin{rem}
	One can arrive to \eqref{main:PointTransformations} either by looking for the best representative for \eqref{eq:NonLinHeath} by utilizing its equivalence algebra to zero out as many arbitrary elements as possible, or by using the Hopf--Cole transformation in combination with a linear transformation since \eqref{eq:NonLinHeath} is very similar to the Burgers equation in potential form and the well established fact in the literature that Eq.~\eqref{eq:Heath} is connected with the Heat equation.
\end{rem}
\begin{rem}
	It is evident that if the function $f$ has the form $f(x,u)= \frac{1}{2}(\exp(\frac{1}{b^2}(a x+u))g(x)+a^2)$, where $g$ an arbitrary real function of $x$, then Eq.~\eqref{eq:NonLinHeath} is linearizable. 
\end{rem}

To classify the class of equations \eqref{main:heatNonLinearSource} we analyzed the classification equation,
\begin{multline*}
	\hat f_\phi \left(\mathcal{F}_1(x,t)+\frac{1}{8} \phi \left(8 \mathcal{F}_4(t)+x \left(x \mathcal{F}_2{}^{\prime \prime }(t)-4 \mathcal{F}_3'(t)\right)\right)\right)+{\mathcal{F}_1}_{,xx}(x,t)-\\
	\frac{1}{8} \hat  f \left(8 \mathcal{F}_4(t)+x \left(x \mathcal{F}_2{}^{\prime \prime }(t)-4 \mathcal{F}_3'(t)\right)\right) +\hat f_x \left(\mathcal{F}_3(t)-\frac{x \mathcal{F}_2'(t)}{2}\right)+\frac{1}{4} \phi \mathcal{F}_2{}^{\prime \prime }(t)-\\
	-\hat f\mathcal{F}_2'(t)-{\mathcal{F}_1}_{,t}(x,t)+\frac{1}{8} \phi \left(x \left(4 \mathcal{F}_3{}^{\prime \prime }(t)-x \mathcal{F}_2{}^{\prime \prime \prime }(t)\right)-8 \mathcal{F}_4'(t)\right)=0,
\end{multline*}
derived from the symmetry condition of \eqref{main:heatNonLinearSource}.

This classification is included in the work of Zhdanov \emph{et al.} \cite{ZhdaLa99}. In their work, the group classification of the heat conductivity equation with a nonlinear source,
\begin{equation*}
	u_t=u_{xx}+F(t,x,u,u_x),
\end{equation*}
was performed taking advantage  of the fact that the abstract Lie algebras of dimensions up to  five are already classified. 

Because of that fact we refrain any further details of the calculations involved and proceed by giving the classification modulo the equivalence group $G^\sim$,
\begin{equation*}
	\tilde x = \delta_4 x +  \delta_3,\ \tilde\tau = \delta _4^2\tau + \delta _0,\ \tilde\phi = \delta _1\phi+\mathcal{F}(x),\ \delta_4\ne0,
\end{equation*}
of the class \eqref{main:heatNonLinearSource}.

In Table \ref{main:table1} the Lie algebra for each form of the function $\hat f$ is given according to the following notations:
{\footnotesize\begin{align*}
	A_1 &= \Span(\partial_\tau),\\
	A_{2,2}^1 &= \Span\left(\partial_{\tau}, e^{\tau } x\partial _x+2e^{\tau}\partial _{\tau} -  e^{\tau}  \left(\frac{1}{4}x^2 +A \right)\phi\partial _\phi\right),\\
	A_{2,2}^2 &= \Span\left(\partial_{\tau}, x\partial _x +2 \tau\partial_{\tau}-B\phi\partial _\phi\right),\\
	A_{2,2}^3 &= \Span\left(\partial_{\tau}, 2e^{2 A\tau}\partial _x-e^{2A\tau} (2 A x + B)\phi\partial _\phi\right),\\
	A_{2,2}^4 &= \Span\left(\partial_{\tau}, e^{-A \tau} \phi\partial _\phi\right),\\
	A_{3,5}^1 &= \Span\left(\partial_{\tau}, e^{2 A\tau}\partial_x-e^{2 A t} \left(A (x+\Delta)\phi+\frac{2 E x^{\frac{2}{B+1}-1} e^{-\frac{1}{2} A (x+\Delta)^2}}{B+1}\right)\partial _\phi,\right.\\
			&\qquad\qquad\quad2 e^{4 A \tau} A (x+\Delta)\partial _x+e^{4 A\tau}\partial _{\tau}-\frac{2 A e^{4 A t}}{B+1} \left( \left(A (B+1) (x+\Delta)^2-2\right)\phi+\right.\\
			&\left.\left.\qquad\qquad 2 \Delta E  x^{\frac{2}{B+1}-1} e^{-\frac{1}{2} A (x+\Delta)^2}\right)\partial _\phi\right),\\
	A_{3,5}^2 &= \Span\left(\partial_{\tau}, e^{\pm\frac{2 \tau}{B+1} }\partial_x \mp \frac{1}{B+1} e^{\pm\frac{2 \tau}{B+1}} x\phi\partial _\phi,\ 2e^{\pm\frac{4 \tau}{B+1} }  x\partial _x\pm e^{\pm\frac{4 \tau}{B+1}}(B+1)\partial _{\tau}\mp\right.\\
			& \left.\frac{ 2 e^{\pm\frac{4 \tau}{B+1}}\left(x^2\mp2\right) \phi}{1+B} \partial _\phi\right),\\
	A_{3,5}^3 &= \Span\Biggl(\partial_{\tau}, \partial_x-\left(\Gamma \phi +\frac{2 \Delta e^{-\Gamma x }x^{\frac{1-B}{1+B}}}{B+1}\right)\partial_\phi,\ (x+2 \Gamma\tau)\partial_x+2\tau\partial_{\tau}+\\
			&\qquad\qquad\ \frac{1}{B+1}\left((2-\Gamma (B+1)(x+2\Gamma\tau ) )\phi-4\Gamma \Delta \tau e^{-\Gamma x} x^{\frac{1-B}{B+1}}\right)\partial_\phi\Biggr),	\\
	A_{3,5}^4 &= \Span\left(\partial_{\tau}, \partial _x-\Gamma\phi\partial _\phi, (x+2 \Gamma\tau )\partial _x+2\tau\partial _{\tau}-\frac{((1+A) \Gamma  (x+2 \Gamma\tau )-2) \phi}{A+1}\partial _\phi\right),\\
	A_{3,5}^5 &= \Span\biggl(\partial_{\tau},  \partial_x-\left(B \phi +\frac{2 e^{-B x}}{ x}\right)\partial_\phi, (x+2B\tau)\partial_x+2\tau\partial_{\tau}-\\
			&\qquad\qquad \left.B\left(\frac{4\tau e^{-B x } }{x}+(x+2 B \tau  ) \phi\right)\partial_\phi\right),\\
	A_{3,5}^6 &= \Span\left(\partial_{\tau},  \partial _x+B\phi\partial _\phi,\ (x-2 B \tau )\partial _x+2 \tau \partial _{\tau}+\left(B(x-2B \tau )\phi- 2 e^{B x }\right)\partial _\phi\right),\\
	A_{3,5}^7 &= \Span\biggl(\partial_{\tau},\ e^{\pm2\tau }\partial _x\pm e^{\pm2\tau }\left(e^{\mp\frac{1}{2}x^2} - \phi\right)x\partial _\phi,\ 2 e^{\pm4\tau} x \partial _x\pm e^{\pm4\tau}\partial_{\tau}+\\
			& \qquad\qquad\  2e^{\pm4\tau}\left(\left(\pm2x^2+3\right)e^{\mp\frac{1}{2} x^2} -  \left(\pm x^2+2\right)\phi\biggr) \partial _\phi\right),\\
	A_{3,5}^8 &= \Span\Biggl(\partial_{\tau}, e^{\pm 2 \tau  } \partial _x- \frac{ e^{\pm 2 t} \left(\pm u x (x+B)+2 e^{\mp\frac{x^2}{2}}\right)}{x+B}\partial _\phi,\\
			& \qquad\qquad\ \pm2e^{\pm4 \tau }x \partial _x+e^{\pm4 \tau }\partial _t+\frac{2  e^{\pm4 t} \left(\pm2 B e^{\mp\frac{x^2}{2}}-u x^2 (x+B)\right)}{x+B}\partial _\phi\Biggr),\\
	A_{3,5}^9 &= \Span\Biggl(\partial_{\tau}, \partial_x-B\phi\partial_\phi, (x+2 B \tau)\partial_x+2\tau\partial_{\tau}-\left(\left(2+B x+2 B^2\tau\right) \phi-B^2e^{-B x} \right)\partial_\phi  \Biggr),\\
	A_{3,5}^{10} &= \Span\left(\partial_{\tau}, e^{\mp2\tau} \partial_x\pm e^{\mp2\tau} x \phi\partial _\phi,\ 2 e^{\mp 4\tau} x\partial _x\mp e^{\mp 4\tau}\partial _t-2 e^{-4 \tau}  \left(2 e^{\pm\frac{1}{2} x^2}\mp x^2\phi\right)\partial _\phi\right),\\
	A_{3,8}^1 	&= \Span\left(\partial_{\tau},  2 \sqrt{B}x \cos(2 \sqrt{B} \tau)\partial _x+2 \sin(2\sqrt{B} \tau )\partial_{\tau}+\right.\\
			&\qquad\qquad\  \left(B \left(\Delta e^{-\frac{1}{2}  \Gamma  x^2}  x^{A+\frac{3}{2}} +x^2\phi\right) \sin(2 \sqrt{B} \tau )+\right.\\
			&\qquad\qquad\  \left.\sqrt{B}  \left(2 \Gamma  \Delta e^{-\frac{1}{2}  \Gamma  x^2}  x^{A+\frac{3}{2}}+ (2 A-1) \phi \right) \cos(2 \sqrt{B} t )\right)\partial _\phi, \\
			&\qquad\qquad\ -2 \sqrt{B} x \sin(2 \sqrt{B} \tau)\partial _x+2 \cos(2\sqrt{B} \tau )\partial_t+\\
			&\qquad\qquad\   \left(B \left(\Delta  e^{-\frac{1}{2}  \Gamma  x^2}x^{A+\frac{3}{2}} + x^2\phi\right) \cos(2 \sqrt{B} \tau )-\right.\\
			&\qquad\quad \left.\left.\sqrt{B}  \left(2 \Gamma  \Delta  e^{-\frac{1}{2}  \Gamma  x^2}x^{A+\frac{3}{2}}+  (2 A-1) \phi\right) \sin(2 \sqrt{B} \tau )\right)\partial _\phi\right),\\	
	A_{3,8}^2 &= \Span\left(\partial_{\tau},  2 \sqrt{\lvert B\rvert} e^{-2 \sqrt{\lvert B\rvert} \tau} x\partial _x-2 e^{-2 \sqrt{\lvert B\rvert} \tau }\partial_{\tau}+\right.\\
			&\qquad\qquad\  e^{-2 \sqrt{\lvert B\rvert}\tau} \sqrt{\lvert B\rvert} \left((\sqrt{\lvert B\rvert}+2 \Gamma) \Delta e^{-\frac{1}{2} \Gamma x^2 } x^{\frac{3}{2}+A} +\left(\sqrt{\lvert B\rvert} x^2+2 A-1\right)\phi\right)\partial _\phi,\\
			&\qquad\qquad\  2 \sqrt{\lvert B\rvert} e^{2 \sqrt{\lvert B\rvert} \tau} x\partial _x+2 e^{2 \sqrt{\lvert B\rvert}\tau }\partial_{\tau}-\\
			&\qquad\qquad\left.\  e^{2  \sqrt{\lvert B\rvert}\tau} \sqrt{\lvert B\rvert} \left((\sqrt{\lvert B\rvert}-2\Gamma) \Delta e^{-\frac{1}{2} \Gamma x^2} x^{\frac{3}{2}+A} + \left(\sqrt{\lvert B\rvert} x^2+1-2 A\right)\phi\right)\partial _\phi\right),	\\
	A_{3,8}^3 &= \Span\biggl(\partial_{\tau},  2 x\partial _x+4 \tau\partial_{\tau}+\left(2 B \Gamma  e^{-\frac{1}{2} B x^2}  x^{\frac{3}{2}+A} +(2 A-1)\phi\right)\partial _\phi, \\
			&\qquad\qquad\ 4 x\tau\partial _x+4\tau^2\partial_{\tau}-\left(\Gamma  e^{-\frac{1}{2} B x^2} (1-4 B \tau ) x^{\frac{3}{2}+A} + \left((2-4A)\tau+x^2\right) \phi\right)\partial _\phi\biggr),\\											
	A_4^1 &= \Span\Biggl(\partial_{\tau},  e^{A\tau } \phi\partial _\phi, 4 e^{\frac{1}{2}  \left(A-\sqrt{A^2-16 B}\right)\tau}\partial _x+\\
			&\qquad\quad\ \ \left(\sqrt{A^2-16 B}-A\right) e^{\frac{1}{2} \left(A-\sqrt{A^2-16 B}\right)\tau}   x \phi\partial _\phi,\ 4 e^{\frac{1}{2}  \left(\sqrt{A^2-16 B}+A\right)\tau}\partial _x-\\
			&\qquad\quad\ \  \left.\ \left(\sqrt{A^2-16 B}+A\right) e^{\frac{1}{2} \left(\sqrt{A^2-16 B}+A\right)\tau}   x \phi\partial _\phi\right),\\
		A_4^2 &= \Span\Biggl(\partial_{\tau},  e^{A\tau } \phi\partial _\phi, 4 e^{\frac{A \tau }{2}} \sin\frac{\sqrt{\lvert A^2-16 B\rvert} \tau }{2}\partial _x-\\
			&\qquad\qquad\ e^{\frac{A \tau }{2}} x \left( \sqrt{\lvert A^2-16 B\rvert}   \cos\frac{ \sqrt{\lvert A^2-16 B\rvert}\tau }{2} + A \sin\frac{ \sqrt{\lvert A^2-16 B\rvert} \tau }{2}  \right)\phi\partial _\phi,\\
			&\qquad\qquad\ 4 e^{\frac{A \tau }{2}} \cos\frac{\sqrt{\lvert A^2-16 B\rvert} \tau }{2}\partial _x+\\
			&\qquad\qquad\ e^{\frac{A \tau }{2}} x\left(\sqrt{\lvert A^2-16 B\rvert}    \sin\frac{ \sqrt{\lvert A^2-16 B\rvert}\tau }{2}-\left. A\cos\frac{ \sqrt{\lvert A^2-16 B\rvert} \tau }{2}  \right)\phi\partial _\phi\right),\\
		A_4^3 &= \Span\left(\partial_{\tau},  e^{\pm\tau } \phi\partial _\phi, 4e^{\frac{\pm\tau}{2}} \partial _x\mp e^{\frac{\pm\tau}{2}} x \phi\partial _\phi,\ 4 e^{\frac{\pm \tau}{2}} \tau \partial _x-e^{\frac{ \tau}{2}} ( 2\pm\tau)x \phi \partial _\phi\right),\\
		A_4^4 &= \Span\left(\partial_{\tau},  e^{A\tau } \phi\partial _\phi, A \partial _x-B \phi\partial _\phi,\ 2 e^{A \tau} \partial _x-e^{A \tau} (A x-2 B \tau)\phi \partial _\phi\right).																				
\end{align*}}%
	\begin{center}
	\footnotesize
	\begin{longtable}{@{}ll@{}}
		\caption[Group classification of Eq.~\eqref{main:heatNonLinearSource}]{Group classification of Eq.~\eqref{main:heatNonLinearSource}}\label{main:table1}\\
		\toprule 
		Lie Algebra    	& $\hat f$\\ \midrule
		\endfirsthead
		\multicolumn{2}{c}{\tablename\ \thetable\ -- \textit{Continued from previous page}} \\ \toprule 
		Lie Algebra    	& $\hat f$\\ \midrule
		\endhead
		\multicolumn{2}{r}{\textit{Continued on next page}} \\
		\bottomrule
		\endfoot
		\bottomrule
		\endlastfoot
			$A_1$      	&  $\forall\hat f$           \\ \addlinespace[0.8em]
			$A_{2,2}^1$	& $-\frac{e^{-\frac{1}{8} x^2}}{4 x^{A+2}} \left(x^2 \left(\frac{1}{4}x^2+2 A-1 \right) \psi +F(\psi)\right),\ \psi=e^{\frac{1}{8} x^2 }x^{A } \phi,\ F^{\prime\prime}\ne0$ \\ \addlinespace[0.8em]
			$A_{2,2}^2$	&  $\frac{1}{x^{A+2}}F(x^A\phi),\ F^{\prime\prime}\ne0$\\ \addlinespace[0.8em]
			$A_{2,2}^3$	&  $-e^{-\frac{1}{2} (A x+B)x}\left(A x (A x + B )\psi+F(\psi)\right),\ \psi=e^{\frac{1}{2}(A x+ B )x} \phi,\ F^{\prime\prime}\ne0$\\ \addlinespace[0.8em]
			$A_{2,2}^4$	&  $-(F(x) +A\log((\Delta+x (\Gamma x +B )) \phi ) )\phi,\ A, B^2+\Gamma^2+\Delta^2\ne0$\\  \addlinespace[0.8em]
			     			& $-e^{-\frac{1}{2} A (x+\Delta)^2} x^{-\frac{2 B}{B+1}}\left( \Gamma\lvert\phi\rvert^{-B}+\frac{A x^2\left(A (B+1) (x+\Delta)^2-B-5\right)}{B+1}\phi-\right.$\\
			$A_{3,5}^1$	& $\left.\qquad\qquad\qquad\qquad\qquad\frac{E \left(A (B+1) x \left(x \left(A (B+1) (x+\Delta)^2-5-B\right)-4 E\right)+2(1- B)\right)}{(B+1)^2} \right)$\\ 
						& $\qquad\quad \psi=e^{\frac{1}{2} A (x+\Delta)^2} x^{-\frac{2}{B+1}} \phi + E,\ A,\Gamma\ne0,\ B\ne0,-1,-2$\\ \addlinespace[0.8em]	
			$A_{3,5}^2$ 	&  $\frac{ 5+B\mp x^2}{(B+1)^2}\phi  - A  e^{\mp\frac{1}{2} x^2} \lvert\phi\rvert ^{-B},\ A \ne0,\ B\ne0,-1,-2$\\ \addlinespace[0.8em]
			\multirow{2}{*}{$A_{3,5}^3$}	& $-e^{-\Gamma x} x^{-\frac{2 B}{1+B}}\left(\frac{\Delta (2 (B+1)-(\Gamma (B+1) x-2)^2)}{(B+1)^2}+ A \lvert\psi\rvert^{-B}+\Gamma^2 x^2 \psi\right),$\\
									& $\qquad\psi=e^{A x} x^{-\frac{2}{1+B}} \phi+\Delta,\ A\ne0,\ B\ne0,-1,-2$\\ \addlinespace[0.8em]
			$A_{3,5}^4$				&  $ -e^{-(A+1)B x}\lvert\phi\rvert^{-A}-B^2\phi,\ A\ne0,-1,-2$\\ \addlinespace[0.8em]
			$A_{3,5}^5$				&  $-B ^2u-\frac{A e^{-e^{B  x}u-B x}}{x^2}-\frac{2 e^{-B x} (2 B  x+1)}{x^2},\ A\ne0$\\\addlinespace[0.8em]
			$A_{3,5}^6$				&  $-e^{B x} \left(A e^\psi+B^2\psi \right),\ \psi=e^{- B x} \phi,\ A\ne0$\\ \addlinespace[0.8em]
		         $A_{3,5}^7$				&  $\mp\frac{e^{\mp\frac{1}{2} x^2}}{4}\left( 4 e^{\pm x^2} u^2 +\left(x^2\mp11\right) \left(x^2\pm1\right) \right)$\\ \addlinespace[0.8em]
			$A_{3,5}^8$				&  $\mp e^{\mp\frac{1}{2} x^2}\left(\frac{A}{(x+B)^2}e^{-\psi }\pm\left(\pm x^2-1\right) \psi+\frac{2(1\mp 2 B(x+B))}{ (x+B)^2}\right),\ \psi= e^{\pm\frac{1}{2} x^2}\phi,\ A\ne0$ \\ \addlinespace[0.8em]
			$A_{3,5}^9$				&  $-e^{B x}  \phi^2 - \frac{ B^4}{4}e^{-Bx}$\\ \addlinespace[0.8em]
			$A_{3,5}^{10}$				&  $ \mp e^{\pm\frac{1}{2} x^2}\left(A e^\psi + (\pm x^2+1)\psi -4\right),\ \psi=e^{\mp\frac{1}{2} x^2} \phi,\ A\ne0$\\ \addlinespace[0.8em]									
									&  $ -e^{-\frac{1}{2} \Gamma x^2 } x^{A-\frac{5}{2}}\psi\left(4A\log\lvert\psi\rvert-\frac{1}{4} x^2 \left(B x^2 +8 A \Gamma\right)\right)+$\\
			$A_{3,8}^1$				&  $\quad\ \ \frac{1}{4} \Delta e^{-\frac{1}{2}\Gamma x^2} x^{A-\frac{5}{2}} \left(3- 8A+ 4( \Gamma x^2-A)^2 \right),\ \psi= e^{\frac{1}{2} \Gamma x^2} x^{\frac{1}{4} (2-4 A)}\phi+\Delta,$\\
									&  $\qquad\quad A\ne0,\ B>0$\\ \addlinespace[0.8em]					
									&  $-e^{-\frac{1}{2} \Gamma x^2 }x^{A-\frac{5}{2}}\psi\left(4A\log\lvert\psi\rvert-\frac{1}{4}x^2 \left(B x^2 +8 A \Gamma\right)\right) +$\\
			$A_{3,8}^2$				&  $\quad\ \ \frac{1}{4} \Delta e^{-\frac{1}{2}\Gamma x^2} x^{A-\frac{5}{2}} \left(3- 8A+ 4( \Gamma x^2-A)^2 \right),\ \psi= e^{\frac{1}{2} \Gamma x^2} x^{\frac{1}{4} (2-4 A)} \phi+\Delta,$\\
									&  $\qquad\quad \ A\ne0,\ B<0$\\ \addlinespace[0.8em]
			\multirow{2}{*}{$A_{3,8}^3$}	&  $2e^{-\frac{1}{2} B x^2 } x^{A-\frac{5}{2}}\psi\left(A B x^2-2A\log\lvert\psi\rvert\right)+\frac{1}{4} \Gamma e^{-\frac{1}{2}B x^2} x^{A-\frac{5}{2}} \left(3- 8A\right)+$\\
									&  $\qquad \Gamma e^{-\frac{1}{2}B x^2} x^{A-\frac{5}{2}} ( B x^2-A)^2,\ \psi= e^{\frac{1}{2} B x^2} x^{\frac{1}{4} (2-4 A)} \phi+\Gamma ,\ A\ne0$\\ \addlinespace[0.8em]								
			$A_4^1$					& $\phi \left(A \log\lvert\phi\rvert + B x^2\right),\ A,B\ne0,\ A^2-16B>0$\\ \addlinespace[0.8em]	
			$A_4^2$			 		& $\phi \left(A \log\lvert\phi\rvert + B x^2\right),\ A,B\ne0,\ A^2-16B<0$\\ \addlinespace[0.8em]
			$A_4^3$			 		& $ \phi \left(\frac{1}{16} x^2\pm \log\lvert\phi\rvert\right)$\\ \addlinespace[0.8em]					
			$A_4^4$					& $\phi  (A \log\lvert\phi\rvert+B x),\ A\ne0$ \\															
	\end{longtable}
	\end{center}
	
In Table \ref{main:table1} we have denoted by $A, B, \Gamma , \Delta, E, ...$  arbitrary real constants, some of them satisfying the stated relations and $\mathcal{F}_1(x)$ is an arbitrary real function.	
\begin{rem}
By investigating the algebraic structure of the symmetry sets we found that, the algebras $A_{3,5}^1$ and $A_{3,5}^2$, the algebras $A_{3,5}^3$ and $A_{3,5}^4$ and the algebras $A_{3,5}^5$ and $A_{3,5}^6$ are connected via an additional equivalence transformation.
\end{rem}

\section{Examples of invariant solutions}

Having obtained the complete group classification for Eq.~\eqref{main:heatNonLinearSource}, and consequently for Eq.~\eqref{eq:Heath}, we can look for invariant solutions under the terminal condition \eqref{cond:terminal} and the barrier option condition \eqref{cond:barrier}: Given a specific Lie algebra from Table~\ref{main:table1}, the appropriate subalgebra, and  the functions $H(t), R(t)$ for the barrier option problem, admitted by each problem are determined using the two required conditions  \eqref{prelim:terminala}, \eqref{prelim:terminalb} and \eqref{prelim:barriera}, \eqref{prelim:barrierb} adapted now to Eq.~\eqref{main:heatNonLinearSource}. Namely,
\begin{equation}\label{main:terminala}
	\left.\mathfrak{X}(\tau-T^\prime)\right\rvert_{\tau=T^\prime}\equiv0,
\end{equation}
\begin{equation}\label{main:terminalb}
	\left.\mathfrak{X}\left(\phi-e^{-\frac{1}{b^2}(ax+1)}\right)\right\rvert_{\tau=T^\prime,\, \phi=e^{-\frac{1}{b^2}(ax+1)}}\equiv0,
\end{equation}
where $T^\prime = -\frac{b^2}{2}T$, and
\begin{equation}\label{main:barriera}
	\left.\mathfrak{X}\left( x-H\left(\tau\right)\right)\right\rvert_{ x=H\left(\tau\right)}\equiv0,
\end{equation}
\begin{equation}\label{main:barrierb}
	\left.\mathfrak{X}\left(\phi-e^{-\frac{1}{b^2}\left(a H\left(\tau\right)+R\left(\tau\right)\right)}\right)\right\rvert_{ x=H\left(\tau\right),\, \phi=e^{-\frac{a}{b^2}H\left(\tau\right)-\frac{1}{b^2}R\left(\tau\right)}}\equiv0.
\end{equation}
Finally, by using that subalgebra similarity solutions are constructed as per usual. 

\subsection{The terminal condition}

Due to the restrictions imposed by the terminal condition all the Lie algebras of dimension less than 4 either they do not yield a suitable subalgebra or the found subalgebra gives trivial similarity solutions (constants or functions depending only on $t$). Hence, we continue by presenting a nontrivial similarity solution for the Lie algebra $A^1_4$.

Let the arbitrary element of the Lie algebra be $\mathfrak X = \mathbf{c}_1 \mathfrak X_1+\mathbf{c}_2 \mathfrak X_2+\mathbf{c}_3 \mathfrak X_3+\mathbf{c}_4 \mathfrak X_4$, where $\mathfrak{X}_i$ are the basis elements of the chosen Lie algebra. Using \eqref{main:terminala} and \eqref{main:terminalb} we obtain the constraints: 
$$
	\mathbf{c}_1=0
$$
and
\begin{multline}
	\frac{1}{4 b^2}e^{-\frac{4+4 a x-2b^2 T^\prime \left(A+\sqrt{A^2-16 B}\right)}{4 b^2}} \left(b^2  \sqrt{A^2-16 B} x(\mathbf{c}_2 -e^{- T^\prime \sqrt{A^2-16 B}}\mathbf{c}_3)-\right.\\
	\left. \left(4 a -b^2 A  x\right) \left(\mathbf{c}_2+e^{-T^\prime \sqrt{A^2-16 B}} \mathbf{c}_3\right)-4 b^2 e^{-\frac{1}{2}  T^\prime \left(\sqrt{A^2-16 B}-A\right)}  \mathbf{c}_4\right)=0.
\end{multline}
From the above two constraints we have the subalgebra defined by 
\begin{align*}
	\mathbf{c}_1&=0,\\
	\mathbf{c}_3&=\frac{e^{ T^\prime \sqrt{A^2-16 B}} \left(A+\sqrt{A^2-16 B}\right) \mathbf{c}_2}{\sqrt{A^2-16 B}-A},\\
	\mathbf{c}_4&=\frac{2 a e^{\frac{1}{2} T^\prime \left(\sqrt{A^2-16 B}-A\right)} \sqrt{A^2-16 B} \mathbf{c}_2}{b^2 \left(A-\sqrt{A^2-16 B}\right)}.
\end{align*}
It follows that the invariant surface condition is  
\begin{multline}\label{main:invariantsurfacecondition1}
	 2 \left(a e^{\frac{1}{2} \left(T^\prime-\tau\right) \left(\sqrt{A^2-16 B}-A\right)} \sqrt{A^2-16 B}+\right.\\
	 \left.2 b^2 \left(e^{ \left(T^\prime-\tau\right) \sqrt{A^2-16 B}}-1\right) B x\right) \phi+ b^2 \left(\left(e^{\left(T^\prime-\tau\right) \sqrt{A^2-16 B}}-1\right) A+\right.\\ \left.\left(1+e^{\left(T^\prime-\tau\right) \sqrt{A^2-16 B}}\right) \sqrt{A^2-16 B}\right) \phi_{x}=0.
\end{multline}
Eq.~\eqref{main:invariantsurfacecondition1} has the solution
\begin{equation}\label{main:similaritysolution1A}
	\phi(x,\tau)= e^{-\frac{2 x \left(a e^{-\frac{1}{2} \left(\tau-T^\prime\right) \left(\sqrt{A^2-16 B}-A\right)} \sqrt{A^2-16 B}+b^2 B\left(e^{ -\left(\tau+T^\prime\right) \sqrt{A^2-16 B}}-1\right) x\right)}{b^2 \left(\left(e^{-\left(\tau-T^\prime\right) \sqrt{A^2-16 B}}-1\right) A+\left(1+e^{ -\left( \tau-T^\prime\right) \sqrt{A^2-16 B}}\right) \sqrt{A^2-16 B}\right)}} \mathcal{F}(\tau),
\end{equation}
where $\mathcal{F}(\tau)$ is an arbritrary real function. Using the above similarity solution Eq.~\eqref{main:heatNonLinearSource} is reduced to the ODE
\begin{multline*}
	\left(2 a^2 e^{\left(\tau-T^\prime\right) \left(A-\sqrt{A^2-16 B}\right)} \left(A^2-16 B\right)+b^4 \left(2 B(\sqrt{A^2-16 B}-A)-\right.\right.\\
	2 e^{-2\left(\tau- T^\prime\right) \sqrt{A^2-16 B}} \sqrt{A^2-16 B} B+A \Delta(A -\sqrt{A^2-16 B}) +\\
	e^{-2\left(\tau-T^\prime\right) \sqrt{A^2-16 B}} A \sqrt{A^2-16 B} \Delta -8 B \Delta +e^{-2\left(\tau-T^\prime\right) \sqrt{A^2-16 B}} \left(A^2 \Delta\right.\\
	 \left.\left.-2 A B-8 B \Delta \right)+4 e^{- \left(\tau-T^\prime\right) \sqrt{A^2-16 B}} \left(A B-4 B \Delta \right)\right)+\\
	 b^4 A \left(\left(1+e^{-2\left(\tau-T^\prime\right) \sqrt{A^2-16 B}}\right) A^2+\left(e^{-2\left(\tau-T^\prime\right) \sqrt{A^2-16 B}}-1\right) A \sqrt{A^2-16 B}-\right.\\
	 \left.  \left.8 \left(1+e^{-\left(\tau- T^\prime\right) \sqrt{A^2-16 B}}\right)^2 B\right) \log\lvert\mathcal{F}\rvert\right) \mathcal{F}+\\
	 b^4 \left(\left(1+e^{-2\left(\tau-T^\prime\right) \sqrt{A^2-16 B}}\right) A^2+\left(e^{-2\left(\tau- T^\prime\right) \sqrt{A^2-16 B}}-1\right) A \sqrt{A^2-16 B}-\right.\\
	 \left.8 \left(1+e^{- \left(\tau-T^\prime\right) \sqrt{A^2-16 B}}\right)^2 B\right) \mathcal{F}^\prime=0.
\end{multline*}
Solving it, by also setting $A=3$ and  $B=1/2$, we obtain the solution
\begin{multline}\label{main:similaritysolution1B}
	\mathcal{F}(\tau) = e^{\frac{e^{3\tau} \left(8 b^4 e^{-4 \tau}+12 a^2 e^{-4 T^\prime}+3 b^4 e^{-\tau-3 T^\prime}-3 e^{-2T^\prime} \left(b^4 e^{-2 \tau}-24 c\right)-2 e^{-\tau-T^\prime} \left(5 b^4 e^{-2 \tau}+72c\right)\right)}{24 b^4 \left(2 e^{-\tau}-e^{-T^\prime}\right)}}\times\\
	 \left(2 e^{-\tau}-e^{-T^\prime}\right)^{-\frac{1}{16} e^{3 (\tau- T^\prime)}},
\end{multline}
where $c$ is the constant of integration. In order the similarity solution found, consisting of Eqs.~\eqref{main:similaritysolution1A} and \eqref{main:similaritysolution1B}, to satisfy the  terminal condition the constant of integration $c$ must be 
$$
	c = \frac{1}{288} e^{-2 T^\prime} \left(48 a^2+96 b^2-8 b^4-6 b^4 T^\prime\right).
$$
Overall, after returning back to the initial Eq.~\eqref{eq:Heath}, we have that the function 
\begin{multline*}
	u(x,t) = \frac{1}{96 b^2 \left(2 e^{\frac{b^2 t}{2}}-e^{\frac{b^2 T}{2}}\right)}e^{-\frac{3 b^2 t}{2}} \left(3 b^6 \left(e^{2 b^2 T}-2 e^{\frac{1}{2} b^2 (t+3 T)}\right) T-\right.\\
	96 a^2 \left(e^{2 b^2 T}-e^{\frac{1}{2} b^2 (t+3 T)}\right)-96 b^2 \left(e^{2 b^2 T}-2 e^{\frac{1}{2} b^2 (t+3 T)}+2 a e^{2 b^2 t} x-\right.\\
	\left.a e^{b^2 (t+T)} x-a e^{\frac{1}{2} b^2 (3 t+T)} x\right)+4 b^4 \left(2 e^{2 b^2 T}+3 e^{b^2 (t+T)}-7 e^{\frac{1}{2} b^2 (t+3 T)}+\right.\\
	\left.4 e^{2 b^2 t} \left(3 x^2-2\right)-2 e^{\frac{1}{2} b^2 (3 t+T)} \left(6 x^2-5\right)\right)-\\
	\left.6 b^4 \left(e^{2 b^2 T}-2 e^{\frac{1}{2} b^2 (t+3 T)}\right) \log\lvert2 e^{\frac{b^2 t}{2}}-e^{\frac{b^2 T}{2}}\rvert\right),
\end{multline*}
is a similarity solution of the PDE
\begin{equation*}
	u_{t}(x,t)=a u_{x}+\frac{1}{2}u_{x}^2-\frac{1}{2}b^2 u_{xx}+\frac{a^2}{2}+\frac{1}{2} b^4 \left(\frac{x^2}{2}-\frac{3 (a x+u)}{b^2}\right).
\end{equation*} 

\subsection{The barrier option condition}
In what follows we present, in detail, a similarity solution using the Lie algebra $A^1_4$ for a specific choice of the respective function $\hat f$. Then, we briefly give examples of similarity solutions for Lie algebras of lower dimension.

Again, let the arbitrary element of the Lie algebra be $\mathfrak X = \mathbf{c}_1 \mathfrak X_1+\mathbf{c}_2 \mathfrak X_2+\mathbf{c}_3 \mathfrak X_3+\mathbf{c}_4 \mathfrak X_4$. Using \eqref{main:barriera} and \eqref{main:barrierb} we obtain the ODE for $H(\tau)$ 
$$
	4 e^{\frac{1}{2}\left(A-\sqrt{A^2-16 B}\right) \tau } \left(\mathbf{c}_3+e^{\sqrt{A^2-16 B}\tau} \mathbf{c}_4\right)-\mathbf{c}_1 H^\prime=0,
$$
with solution
\begin{multline}\label{main:solutionOfH}
	H(\tau)=\frac{e^{\frac{1}{2} \left(A-\sqrt{A^2-16 B}\right)\tau} }{2 B \mathbf{c}_1}\left(\sqrt{A^2-16 B}( \mathbf{c}_3-e^{ \sqrt{A^2-16 B}\tau} \mathbf{c}_4)+\right.\\
		 \left.A \left(\mathbf{c}_3+e^{\sqrt{A^2-16 B}\tau} \mathbf{c}_4\right)+2 e^{-\frac{1}{2} \left(A-\sqrt{A^2-16 B}\right)\tau} B \mathbf{c}_1 \mathbf{c}_{5}\right),\, \mathbf{c}_{1}\ne0,
\end{multline}
where $\mathbf{c}_{5}$ is the constant of integration, and an ODE for $R(\tau)$
\begin{multline*}
	 4 a  B \mathbf{c}_1 e^{\frac{1}{2} \left(A-\sqrt{A^2-16 B}\right)\tau}\left(\mathbf{c}_3+e^{\sqrt{A^2-16 B}\tau} \mathbf{c}_4\right)+\\
	 b^2 \left(B \left( \mathbf{c}_1 \mathbf{c}_3  \mathbf{c}_5 \left(\sqrt{A^2-16 B}-A\right) e^{\frac{1}{2}  \left(A-\sqrt{A^2-16 B}\right)\tau}-\right.\right.\\
	 \mathbf{c}_1 \mathbf{c}_4 \mathbf{c}_5 \left(\sqrt{A^2-16 B}+A\right) e^{\frac{1}{2}  \left(\sqrt{A^2-16 B}+A\right)\tau}-8 \mathbf{c}_3{}^2 e^{ \left(A-\sqrt{A^2-16 B}\right)\tau}-\\
	  \left.\left.8 \mathbf{c}_4{}^2 e^{ \left(\sqrt{A^2-16 B}+A\right)\tau}+\left(\mathbf{c}_1 \mathbf{c}_2+16 \mathbf{c}_3 \mathbf{c}_4\right) e^{A \tau}\right)-2 A^2 \mathbf{c}_3 \mathbf{c}_4 e^{A \tau}\right) +  B \mathbf{c}_1^2R^\prime(t)=0,
\end{multline*}
with solution
\begin{multline}\label{main:solutionOfR}
	R(\tau)=\frac{1}{B \mathbf{c}_1^2}\left(\frac{1}{2} e^{\frac{1}{2}  \left(A-\sqrt{A^2-16 B}\right)\tau} \left(b^2 \left(4 B \mathbf{c}_1 \left(\mathbf{c}_4\mathbf{c}_5 e^{\sqrt{A^2-16 B}\tau}+\mathbf{c}_3\right)+\right.\right.\right.\\
			e^{\frac{1}{2}  \left(A-\sqrt{A^2-16 B}\right)\tau} \left(A \left(\mathbf{c}_4{}^2 e^{2 \sqrt{A^2-16 B}\tau}+\mathbf{c}_3{}^2\right)+\right.\\
			\left.\left. \sqrt{A^2-16 B} \left(\mathbf{c}_3{}^2-\mathbf{c}_4{}^2 e^{2  \sqrt{A^2-16 B}\tau}\right)\right)\right)-a \mathbf{c}_1 \left(A \left(\mathbf{c}_4 e^{ \sqrt{A^2-16 B}\tau}+\mathbf{c}_3\right)+\right.\\
			\left.\left. \sqrt{A^2-16 B} \left(\mathbf{c}_3-\mathbf{c}_4 e^{\sqrt{A^2-16 B}\tau}\right)\right)\right)-\frac{b^2 B \left(\mathbf{c}_1 \mathbf{c}_2+16 \mathbf{c}_3 \mathbf{c}_4\right) e^{A \tau}}{A}+\\
			\left.2 A b^2 \mathbf{c}_3 \mathbf{c}_4 e^{A \tau}\right)+\mathbf{c}_6,
\end{multline}
where $\mathbf{c}_{6}$ is the constant of integration.

The solution  $H$ found above includes the barrier function usually used in literature, \eqref{H0}. Indeed, by setting
\begin{align*}
	\mathbf{c}_1	&=-\frac{2 b^2 e^{\alpha T}\mathbf{c}_3 }{\alpha\beta K},\ \mathbf{c}_4=\mathbf{c}_5=0
	\intertext{and}
	B 			&=-\frac{b^2 \alpha  A+2 \alpha ^2}{2 b^4},
\end{align*}
the functions $H$ and $R$ became respectively
\begin{align*}
	H(\tau) 	&=K \beta e^{-\frac{2 \alpha(\tau + b^2 T/2)}{b^2}} \\
	\intertext{and}	
	R(\tau)	&=\frac{1}{2} \beta  K \left(\frac{\alpha  \mathbf{c}_2 e^{A \tau-\alpha  T}}{A \mathbf{c}_3}-e^{-2 \alpha  \left(\frac{2 \tau}{b^2}+T\right)} \left(2 a e^{\alpha  \left(\frac{2 \tau}{b^2}+T\right)}+\alpha  \beta  K\right)\right)+\mathbf{c}_6,
\end{align*}
where we have assumed that $A b^2>-4\alpha$. Furthermore, by setting 
\begin{align*}
	\mathbf{c}_2 	&=\mathbf{c}_6=0
\end{align*}
$R$ simplifies to
$$
	R(\tau) =-\frac{1}{2} \beta  K e^{-2 \alpha  \left(\frac{2 \tau}{b^2}+T\right)} \left(2 a e^{\alpha  \left(\frac{2 \tau}{b^2}+T\right)}+\alpha  \beta  K\right).
$$

It follows that the invariant surface condition is  now
\begin{equation}\label{main:invariantsurfacecondition2}
	b^4 e^{\frac{\alpha(2\tau +b^2 T)}{b^2}} \phi_{\tau}+2 \alpha  \beta  K \left(\alpha x \phi-b^2 \phi_x\right)=0.
\end{equation}
Eq.~\eqref{main:invariantsurfacecondition2} has the solution
\begin{equation}\label{main:similaritysolution2A}
	\phi(x,\tau)= e^{\frac{ \alpha x^2}{2 b^2}} \mathcal{F}\left(\frac{b^2 e^{-\frac{2  \alpha\tau }{b^2}-\alpha T} \left(e^{\frac{2 \alpha\tau }{b^2}+\alpha T} x-K \beta \right)}{2 \alpha  \beta K}\right).
\end{equation}
Using the above similarity solution Eq.~\eqref{main:heatNonLinearSource} is reduced to the ODE
\begin{equation*}
4 \alpha ^2 \beta ^2 K^2 \left(2 \alpha  \zeta  \mathcal{F}^\prime+\mathcal{F} \left(\alpha +A b^2 \log\lvert\mathcal{F}\rvert+b^2 \Delta \right)\right)+b^6 \log\lvert\mathcal{F}\rvert^{\prime \prime }=0,
\end{equation*}
where $\zeta=\frac{b^2 e^{-\frac{2 \alpha  t}{b^2}-\alpha T} \left(x e^{\frac{2 \alpha  t}{b^2}+\alpha  T}-\beta  K\right)}{2 \alpha  \beta  K}$. A special solution of the above ODE is the following
\begin{equation}\label{main:similaritysolution2B}
	\mathcal{F}(\zeta) = \exp\left(\frac{2 \alpha +A b^2-2 b^2 \Delta }{2 A b^2}-\frac{\alpha ^2 \beta ^2 K^2 \left(4 \alpha +A b^2\right)}{b^6} \zeta ^2\right).
\end{equation}
In addition, for the above solution to satisfy the condition $\phi=e^{-\frac{a H\left(\tau\right)+R\left(\tau\right)}{b^2}}$ the constant $\Delta$ must be set equal to  $\frac{A}{2}+\frac{\alpha }{b^2}$.

Finally, after returning back to the initial Eq.~\eqref{eq:Heath}, we have that the function 
\begin{equation*}
	u(x,t) =\frac{1}{4} \left(\left(4 \alpha +A b^2\right) \left(x-\beta  K e^{\alpha  (t-T)}\right)^2-2 \alpha  x^2\right)-a x,
\end{equation*}
is a similarity solution of the PDE
\begin{multline*}
	u_{t}(x,t)=a u_{x}+\frac{1}{2}u_{x}^2-\frac{1}{2}b^2 u_{xx}+\frac{1}{2}a^2-\\
	\frac{1}{2}  A b^2 (a x+u)+\frac{1}{4}\left(2 \alpha +A b^2\right) \left(b^2-\alpha  x^2\right).
\end{multline*} 

By using now the Lie algebra $A^1_{2,2}$ we can obtain the similarity solution 
\begin{equation*}
	u(x,t)=b^2 \left(-\log \left(e^{-\frac{x^2}{8}} \sqrt{x} \log \left(8 \sec ^2\left(\frac{1}{4} \left(b^2 t-2 c_3\right)+\log x\right)\right)\right)\right)-a x
\end{equation*}
for the PDE
\begin{multline*}
	u_{t}(x,t)=a u_{x}+\frac{1}{2}u_{x}^2-\frac{1}{2}b^2 u_{xx}+\\
	\frac{1}{32} \left(16 a^2+b^4 \left(8-\frac{4 \sqrt{x} e^{\frac{e^{\frac{x^2}{8}-\frac{a x+u}{b^2}}}{\sqrt{x}}+\frac{a x+u}{b^2}-\frac{x^2}{8}}+x^4-4}{x^2}\right)\right).
\end{multline*} 

And finally by using the Lie algebra $A^9_{3,5}$ we can obtain the similarity solution 
\begin{equation*}
	u(x,t)=b^2 \left(-\log \left(\frac{3}{2} e^{-3 x} \left(\frac{4}{\left(3 b^2 t-6 c_1+x\right){}^2}+3\right)\right)\right)-a x
\end{equation*}
for the PDE
\begin{multline*}
	u_{t}(x,t)=a u_{x}+\frac{1}{2}u_{x}^2-\frac{1}{2}b^2 u_{xx}+\\
	\frac{1}{8} \left(4 a^2-b^4 \left(77 \sinh \left(\frac{a x-3 b^2 x+u}{b^2}\right)+85 \cosh \left(\frac{a x-3 b^2 x+u}{b^2}\right)\right)\right).
\end{multline*} 

\section{Conclusion}

In the present paper a generalization of the Heath equation \eqref{eq:Heath}  was proposed and studied under the view of the modern group analysis. To that end, we harnessed the advantage of being able to connect it with the heat equation with nonlinear course, a well known and studied equation. This fact substantially simplifies the task of classifying it and obtaining its point symmetries.  

Through this classification interesting cases, from the point of view of symmetries, arise. Nonlinear equations in general have few or no symmetries so cases that augment the set of symmetries at disposal are like an oasis in the desert. After all, it is evident in the related  literature that a dynamical system possessing an ample number of symmetries is more probable to relate with a physical system or model a more realistic process. Furthermore, in the case that we wish to study a boundary problem, because of the fact that not all of the symmetries admit the boundary conditions, some of the symmetries will be excluded. Hence the bigger the set of symmetries the bigger the probability that some will survive the scrutiny of the boundary conditions and give an invariant solution for the problem in its entirety. 

This is evident for the terminal condition where due to the restrictions imposed by it only four dimensional algebras are able to yield nontrivial solutions.

Things are different when the barrier option is considered, because of the two functions $H(t),R(t)$ a broader range of cases can yield interesting solutions. It is worth mentioning at this point that, as can be seen by the second example, the barrier function $H$ usually used in the related literature is admitted by the symmetries. This fact  reinforces further our belief that symmetries methods can be a valuable tool in investigating this kind of financial problems.

Indeed, the insight provided through the above symmetry analysis might prove practical to anyone looking for a more realistic economic model without departing from the reasoning behind the Heath model. Moreover,  when one studies more exotic kinds of options, e.g. options that have gained ground in the Asian markets which in turn play an ever increasing role in the world market. The nonlinear variants of the traditional model given here might turn the table in that respect. We leave to the interested reader the possible economical interpretation and use of the obtained results.

\section*{Acknowledgements} 

 We would like to thank the reviewers for their useful remarks that helped the manuscript to reach its present form. We would also like to thank FAPESP for the PostDoc grant (Proc. \#2011/05855-9)  giving S. Dimas the opportunity to visit IMECC--UNICAMP where this work was carried out.

\bibliographystyle{model1-num-names}
\bibliography{Bibliography}

\begin{thebibliography}{39}
\expandafter\ifx\csname natexlab\endcsname\relax\def\natexlab#1{#1}\fi
\providecommand{\bibinfo}[2]{#2}
\ifx\xfnm\relax \def\xfnm[#1]{\unskip,\space#1}\fi
\bibitem[{Black and Scholes(1973)}]{BlaScho73}
\bibinfo{author}{F.~Black}, \bibinfo{author}{M.~Scholes},
\newblock \bibinfo{title}{The pricing of options and corporate liabilities},
\newblock \bibinfo{journal}{J. Polit. Econ.} \bibinfo{volume}{81}
  (\bibinfo{year}{1973}) \bibinfo{pages}{637--659}.
\bibitem[{Longstaff(1989)}]{Lo89}
\bibinfo{author}{F.~A. Longstaff},
\newblock \bibinfo{title}{A nonlinear general equilibrium model of the term
  structure of interest rates},
\newblock \bibinfo{journal}{J. Financ. Econ.} \bibinfo{volume}{23}
  (\bibinfo{year}{1989}) \bibinfo{pages}{195--224}.
\bibitem[{Vasicek(1977)}]{Va77}
\bibinfo{author}{O.~Vasicek},
\newblock \bibinfo{title}{An equilibrium characterization of the term
  structure},
\newblock \bibinfo{journal}{J. Financ. Econ.} \bibinfo{volume}{5}
  (\bibinfo{year}{1977}) \bibinfo{pages}{177--188}.
\bibitem[{Cox et~al.(1985)Cox, Ingersoll, and Ross}]{CoIngeRo85}
\bibinfo{author}{J.~C. Cox}, \bibinfo{author}{J.~E. Ingersoll},
  \bibinfo{author}{S.~A. Ross},
\newblock \bibinfo{title}{An intertemporal general equilibrium model of asset
  prices},
\newblock \bibinfo{journal}{Econometrica} \bibinfo{volume}{53}
  (\bibinfo{year}{1985}) \bibinfo{pages}{363--384}.
\bibitem[{Heath et~al.(2001)Heath, Platin, and Schweizer}]{HeaPlaSchwei2k1}
\bibinfo{author}{D.~Heath}, \bibinfo{author}{E.~Platin},
  \bibinfo{author}{M.~Schweizer},
\newblock \bibinfo{title}{Numerical comparison of local risk-minimisation and
  mean-variance hedging},
\newblock in: \bibinfo{editor}{E.~Jouini}, \bibinfo{editor}{C.~Jajusa},
  \bibinfo{editor}{M.~Murek} (Eds.), \bibinfo{booktitle}{Option Pricing,
  Interest Rates and Risk Management}, \bibinfo{publisher}{CUP, Cambridge},
  \bibinfo{year}{2001}, pp. \bibinfo{pages}{509--537}.
\bibitem[{Gazizov and Ibragimov(1998)}]{GazIbr98}
\bibinfo{author}{R.~K. Gazizov}, \bibinfo{author}{N.~H. Ibragimov},
\newblock \bibinfo{title}{Lie symmetry analysis of differential equations in
  finance},
\newblock \bibinfo{journal}{Nonlinear Dynam.} \bibinfo{volume}{17}
  (\bibinfo{year}{1998}) \bibinfo{pages}{387--407}.
\bibitem[{Dimas et~al.(2009)Dimas, Andriopoulos, Tsoubelis, and
  Leach}]{DiAndTsLe2k9}
\bibinfo{author}{S.~Dimas}, \bibinfo{author}{K.~Andriopoulos},
  \bibinfo{author}{D.~Tsoubelis}, \bibinfo{author}{P.~G.~L. Leach},
\newblock \bibinfo{title}{Complete specification of some partial differential
  equations that arise in financial mathematics},
\newblock \bibinfo{journal}{J. Nonlinear Math. Phys.} \bibinfo{volume}{16, s-1}
  (\bibinfo{year}{2009}) \bibinfo{pages}{73--92}.
\bibitem[{Sinkala et~al.(2008)Sinkala, Leach, and O'Hara}]{SiLeaHa2k8}
\bibinfo{author}{O.~Sinkala}, \bibinfo{author}{P.~Leach},
  \bibinfo{author}{J.~O'Hara},
\newblock \bibinfo{title}{Invariance properties of a general bond-pricing
  equation},
\newblock \bibinfo{journal}{J. Diff. Eq.} \bibinfo{volume}{244}
  (\bibinfo{year}{2008}) \bibinfo{pages}{2820--2835}.
\bibitem[{Bozhkov and Dimas(2014)}]{BoDi2k13x}
\bibinfo{author}{Y.~Bozhkov}, \bibinfo{author}{S.~Dimas}, \bibinfo{title}{Group
  analysis of a semi-linear general bond-pricing equation},
  \bibinfo{year}{2014}. \bibinfo{note}{Submitted}.
\bibitem[{Naicker et~al.(2005)Naicker, Andriopoulos, and
  Leach}]{NaiAndrioLea2k5}
\bibinfo{author}{V.~Naicker}, \bibinfo{author}{K.~Andriopoulos},
  \bibinfo{author}{P.~G.~L. Leach},
\newblock \bibinfo{title}{Symmetry {R}eductions of a
  {H}amilton--{J}acobi--{B}ellman {E}quation {A}rising in {F}inancial
  {M}athematics},
\newblock \bibinfo{journal}{J. Nonlinear Math. Phys.} \bibinfo{volume}{12}
  (\bibinfo{year}{2005}) \bibinfo{pages}{268--283}.
\bibitem[{Bluman and Kumei(1989)}]{BluKu89}
\bibinfo{author}{G.~W. Bluman}, \bibinfo{author}{S.~Kumei},
  \bibinfo{title}{Symmetries and differential equations},
  \bibinfo{publisher}{Springer}, \bibinfo{address}{New York},
  \bibinfo{year}{1989}.
\bibitem[{Ovsiannikov(1982)}]{Ovsi82}
\bibinfo{author}{L.~Ovsiannikov}, \bibinfo{title}{Group Analysis of
  Differential Equations}, \bibinfo{publisher}{Academic Press},
  \bibinfo{edition}{1$^{st}$} edition, \bibinfo{year}{1982}. \bibinfo{note}{432
  pages}.
\bibitem[{Ibragimov(2009)}]{Ibra2k9}
\bibinfo{author}{N.~H. Ibragimov},
\newblock \bibinfo{title}{Equivalence groups and invariants of linear and
  nonlinear equations},
\newblock \bibinfo{journal}{Archives of ALGA} \bibinfo{volume}{4}
  (\bibinfo{year}{2009}) \bibinfo{pages}{41--100}.
\bibitem[{Popovych and Eshraghi(2005)}]{PoEshra2k5}
\bibinfo{author}{R.~O. Popovych}, \bibinfo{author}{H.~Eshraghi},
\newblock \bibinfo{title}{Admissible point transformations of nonlinear
  {S}chr\"odinger equations},
\newblock in: \bibinfo{editor}{N.~Ibragimov}, \bibinfo{editor}{C.~Sophocleous},
  \bibinfo{editor}{P.~Damianou} (Eds.), \bibinfo{booktitle}{Proceedings of the
  10$^{th}$ International Conference in {MO}dern {GR}oup {AN}alysis},
  \bibinfo{year}{2005}, pp. \bibinfo{pages}{167--174}.
\bibitem[{Romano and Torrisi(1999)}]{RoTo99}
\bibinfo{author}{V.~Romano}, \bibinfo{author}{M.~Torrisi},
\newblock \bibinfo{title}{Application of weak equivalence transformations to a
  group analysis of a drift-diffusion model},
\newblock \bibinfo{journal}{J. Phys. A: Math. Gen.} \bibinfo{volume}{32}
  (\bibinfo{year}{1999}) \bibinfo{pages}{7953}.
\bibitem[{Bozhkov and Dimas(2014)}]{BoDi2k13b}
\bibinfo{author}{Y.~Bozhkov}, \bibinfo{author}{S.~Dimas},
\newblock \bibinfo{title}{Group classification of a generalized
  {B}lack--{S}choles--{M}erton equation},
\newblock \bibinfo{journal}{Commun. Nonlinear Sci. Numer. Simul.}
  \bibinfo{volume}{19} (\bibinfo{year}{2014}) \bibinfo{pages}{2200--2211}.
\bibitem[{Cherniha et~al.(2008)Cherniha, Serov, and Rassokha}]{CheSeRa2k8}
\bibinfo{author}{R.~Cherniha}, \bibinfo{author}{M.~Serov},
  \bibinfo{author}{I.~Rassokha},
\newblock \bibinfo{title}{Lie symmetries and form-preserving transformations of
  reaction--diffusion--convection equations},
\newblock \bibinfo{journal}{J. Math. Anal. Appl.} \bibinfo{volume}{342}
  (\bibinfo{year}{2008}) \bibinfo{pages}{1363--1379}.
\bibitem[{Cardoso-Bilho et~al.(2011)Cardoso-Bilho, Bilho, and
  Popovych}]{CaBiPo2k11}
\bibinfo{author}{E.~D.~S. Cardoso-Bilho}, \bibinfo{author}{A.~Bilho},
  \bibinfo{author}{R.~O. Popovych},
\newblock \bibinfo{title}{Enhanced preliminary group classification of a class
  of generalized diffusion equations},
\newblock \bibinfo{journal}{Commun. Nonlinear Sci. Numer. Simul.}
  \bibinfo{volume}{16} (\bibinfo{year}{2011}) \bibinfo{pages}{3622--3638}.
\bibitem[{Ivanova et~al.(2010)Ivanova, Popovych, and Sophocleous}]{IvaPoSo2k10}
\bibinfo{author}{N.~M. Ivanova}, \bibinfo{author}{R.~O. Popovych},
  \bibinfo{author}{C.~Sophocleous},
\newblock \bibinfo{title}{Group analysis of variable coefficient
  diffusion-convection equations. {I}. {E}nhanced group classification},
\newblock \bibinfo{journal}{Lobachevskii J. Math.} \bibinfo{volume}{31}
  (\bibinfo{year}{2010}) \bibinfo{pages}{100--122}.
\bibitem[{Vaneeva et~al.(2009)Vaneeva, Popovych, and Sophocleous}]{VaPoSo2k9}
\bibinfo{author}{O.~O. Vaneeva}, \bibinfo{author}{R.~O. Popovych},
  \bibinfo{author}{C.~Sophocleous},
\newblock \bibinfo{title}{Enhanced group analysis and exact solutions of
  variable coefficient semilinear diffusion equations with a power source},
\newblock \bibinfo{journal}{Acta Appl. Math.} \bibinfo{volume}{106}
  (\bibinfo{year}{2009}) \bibinfo{pages}{1--46}.
\bibitem[{Black and Scholes(1972)}]{BlaScho72}
\bibinfo{author}{F.~Black}, \bibinfo{author}{M.~Scholes},
\newblock \bibinfo{title}{The valuation of option contracts and a test of
  market efficiency},
\newblock \bibinfo{journal}{J. Finance} \bibinfo{volume}{27}
  (\bibinfo{year}{1972}) \bibinfo{pages}{399--417}.
\bibitem[{Merton(1974)}]{Me74}
\bibinfo{author}{R.~C. Merton},
\newblock \bibinfo{title}{On the pricing of corporate debt: The risk structure
  of interest rates},
\newblock \bibinfo{journal}{J. Finance} \bibinfo{volume}{29}
  (\bibinfo{year}{1974}) \bibinfo{pages}{449--470}.
\bibitem[{Ugur(2008)}]{Ugur2k8}
\bibinfo{author}{O.~Ugur}, \bibinfo{title}{Introduction to computational
  finance}, \bibinfo{publisher}{Imperial College Press and World Scientific},
  \bibinfo{year}{2008}.
\bibitem[{O'Hara(2011)}]{ha2k11}
\bibinfo{author}{J.~O'Hara}, \bibinfo{title}{Lecture notes on exotic options},
  \bibinfo{year}{2011}. \bibinfo{note}{Http://courses.essex.ac.uk/cf/cf966/}.
\bibitem[{Kwok(2008)}]{Kwo2k8}
\bibinfo{author}{Y.~K. Kwok}, \bibinfo{title}{Mathematical Models of Financial
  Derivatives}, \bibinfo{publisher}{Springer}, \bibinfo{edition}{2nd} edition,
  \bibinfo{year}{2008}.
\bibitem[{O'Hara et~al.(2013)O'Hara, Sophocleous, and Leach}]{HaSoLea2k13}
\bibinfo{author}{J.~O'Hara}, \bibinfo{author}{C.~Sophocleous},
  \bibinfo{author}{P.~G.~L. Leach},
\newblock \bibinfo{title}{Symmetry analysis of a model for the exercise of a
  barrier option},
\newblock \bibinfo{journal}{Commun. Nonlinear Sci. Numer. Simul.}
  \bibinfo{volume}{18} (\bibinfo{year}{2013}) \bibinfo{pages}{2367--2373}.
\bibitem[{Olver(2000)}]{Olver2k}
\bibinfo{author}{P.~J. Olver}, \bibinfo{title}{Applications of Lie Groups to
  Differential Equations}, volume \bibinfo{volume}{107} of
  \textit{\bibinfo{series}{Graduate Texts in Mathematics}},
  \bibinfo{publisher}{Springer}, \bibinfo{address}{New York},
  \bibinfo{edition}{2$^{nd}$} edition, \bibinfo{year}{2000}.
\bibitem[{Head(1993)}]{Head93}
\bibinfo{author}{A.~K. Head},
\newblock \bibinfo{title}{Lie, a pc program for {L}ie analysis of differential
  equations},
\newblock \bibinfo{journal}{Comput. Phys. Comm.} \bibinfo{volume}{77}
  (\bibinfo{year}{1993}) \bibinfo{pages}{241--248}.
\bibitem[{Nucci(1996)}]{Nucci1}
\bibinfo{author}{M.~Nucci},
\newblock \bibinfo{title}{Interactive {REDUCE} programs for calculating {L}ie
  point, non-classical, {L}ie-{B}{\"a}cklund, and approximate symmetries of
  differential equations: manual and floppy disk},
\newblock in: \bibinfo{editor}{N.~Ibragimov} (Ed.), \bibinfo{booktitle}{CRC
  Handbook of Lie Group Analysis of Differential Equations. Vol. III:New
  Trends}, \bibinfo{publisher}{CRC Press}, \bibinfo{year}{1996}, pp.
  \bibinfo{pages}{415--481}.
\bibitem[{Nucci(1992)}]{Nucci2}
\bibinfo{author}{M.~Nucci},
\newblock \bibinfo{title}{Interactive {REDUCE} programs for calculating
  classical, nonclassical and {L}ie-{B}\"acklund symmetries for differential
  equations},
\newblock in: \bibinfo{editor}{W.~Ames}, \bibinfo{editor}{P.~Van~der Houwen}
  (Eds.), \bibinfo{booktitle}{Computational and Applied Mathematics II.
  Differential Equations}, \bibinfo{publisher}{Elsevier}, \bibinfo{year}{1992},
  pp. \bibinfo{pages}{345--350}.
\bibitem[{Baumann(2000)}]{Baumann2k}
\bibinfo{author}{G.~Baumann}, \bibinfo{title}{Symmetry Analysis of Differential
  Equations with Mathematica}, \bibinfo{publisher}{Telos/Springer},
  \bibinfo{address}{New York}, \bibinfo{year}{2000}.
\bibitem[{Reasearch{, Inc.}(2010)}]{Wolfram2k10}
\bibinfo{author}{W.~Reasearch{, Inc.}}, \bibinfo{title}{Mathematica Edition:
  Version 8.0}, \bibinfo{publisher}{Wolfram Reasearch, Inc.},
  \bibinfo{address}{Champaign, Illinois}, \bibinfo{year}{2010}.
\bibitem[{Dimas(2008)}]{Dimas2k8}
\bibinfo{author}{S.~Dimas}, \bibinfo{title}{Partial differential equations,
  algebraic computing and nonlinear systems}, \bibinfo{type}{Ph.{D}. {T}hesis},
  University of Patras, \bibinfo{address}{Patras, Greece},
  \bibinfo{year}{2008}.
\bibitem[{Dimas and Tsoubelis(2005)}]{DiTs2k5a}
\bibinfo{author}{S.~Dimas}, \bibinfo{author}{D.~Tsoubelis},
\newblock \bibinfo{title}{{SYM}: A new symmetry-finding package for
  {M}athematica},
\newblock in: \bibinfo{editor}{N.~Ibragimov}, \bibinfo{editor}{C.~Sophocleous},
  \bibinfo{editor}{P.~Damianou} (Eds.), \bibinfo{booktitle}{The 10$^{th}$
  International Conference in {MO}dern {GR}oup {AN}alysis},
  \bibinfo{publisher}{University of Cyprus}, \bibinfo{address}{Nicosia},
  \bibinfo{year}{2005}, pp. \bibinfo{pages}{64--70}.
\bibitem[{Dimas and Tsoubelis(2006)}]{DiTs2k6}
\bibinfo{author}{S.~Dimas}, \bibinfo{author}{D.~Tsoubelis},
\newblock \bibinfo{title}{A new {M}athematica-based program for solving
  overdetermined systems of {PDEs}},
\newblock in: \bibinfo{editor}{Y.~Papegay} (Ed.), \bibinfo{booktitle}{Applied
  {Mathematica}, {E}lectronic {P}roceedings of the {E}ighth {I}nternational
  {M}athematica {S}ymposium (IMS'06)}, \bibinfo{publisher}{France: INRIA},
  \bibinfo{address}{Avignon, France}, \bibinfo{year}{2006}. \bibinfo{note}{ISBN
  2-7261-1289-7}.
\bibitem[{Ibragimov(1985)}]{Ibra85b}
\bibinfo{author}{N.~H. Ibragimov}, \bibinfo{title}{Transformation Groups
  Applied to Mathematical Physics}, Mathematics and its Applications,
  \bibinfo{publisher}{Springer}, \bibinfo{edition}{1$^{st}$} edition,
  \bibinfo{year}{1985}.
\bibitem[{Hydon(2000)}]{Hydon2k}
\bibinfo{author}{P.~E. Hydon}, \bibinfo{title}{Symmetry Methods for
  Differential Equations}, Cambridge Texts in Applied Mathematics,
  \bibinfo{publisher}{Cambridge University Press},
  \bibinfo{address}{Cambridge}, \bibinfo{edition}{1$^{st}$} edition,
  \bibinfo{year}{2000}.
\bibitem[{Stephani(1990)}]{Ste90}
\bibinfo{author}{H.~Stephani}, \bibinfo{title}{Differential Equations: Their
  Solution Using Symmetries}, \bibinfo{publisher}{Cambridge University Press},
  \bibinfo{address}{Cambridge}, \bibinfo{edition}{1$^{st}$} edition,
  \bibinfo{year}{1990}. \bibinfo{note}{Editor: MacCallum, Malcolm}.
\bibitem[{Zhdanov and Lahno(1999)}]{ZhdaLa99}
\bibinfo{author}{R.~Z. Zhdanov}, \bibinfo{author}{V.~I. Lahno},
\newblock \bibinfo{title}{Group classification of heat conductivity equations
  with a nonlinear source.},
\newblock \bibinfo{journal}{J. Phys. A: Math. Gen.} \bibinfo{volume}{32}
  (\bibinfo{year}{1999}) \bibinfo{pages}{7405--7418}.

\end{thebibliography}

\end{document}